\documentclass{amsart}
\usepackage{amsfonts}

\newtheorem{theorem}{Theorem}
\theoremstyle{plain}

\newtheorem{corollary}{Corollary}

\newtheorem{lemma}{Lemma}

\newtheorem{proposition}{Proposition}
\newtheorem{remark}{Remark}

\numberwithin{equation}{section}
\input{tcilatex}

\begin{document}
\title[Bessel's Inequality in Inner Product Spaces]{A Counterpart of
Bessel's Inequality in Inner Product Spaces and Some Gr\"{u}ss Type Related
Results}
\author{S.S. Dragomir}
\address{School of Computer Science and Mathematics\\
Victoria University of Technology\\
PO Box 14428, MCMC 8001\\
Victoria, Australia.}
\email{sever.dragomir@vu.edu.au}
\urladdr{http://rgmia.vu.edu.au/SSDragomirWeb.html}
\date{May 16, 2003}

\begin{abstract}
A counterpart of the famous Bessel's inequality for orthornormal families in
real or complex inner product spaces is given. Applications for some Gr\"{u}%
ss type inequalities are also provided.
\end{abstract}

\maketitle

\section{Introduction\label{s1}}

In \cite{SSD1}, the author has proved the following Gr\"{u}ss type
inequality in real or complex inner product spaces.

\begin{theorem}
\label{t1}Let $\left( H,\left\langle \cdot ,\cdot \right\rangle \right) $ be
an inner product space over $\mathbb{K}$ $\left( \mathbb{K}=\mathbb{R},%
\mathbb{C}\right) $ and $e\in H,$ $\left\Vert e\right\Vert =1.$ If $\phi
,\Phi ,\gamma ,\Gamma $ are real or complex numbers and $x,y$ are vectors in 
$H$ such that the conditions%
\begin{equation}
\func{Re}\left\langle \Phi e-x,x-\phi e\right\rangle \geq 0\text{ \ and \ }%
\func{Re}\left\langle \Gamma e-y,y-\gamma e\right\rangle \geq 0  \label{1.1}
\end{equation}%
hold, then we have the inequality%
\begin{equation}
\left\vert \left\langle x,y\right\rangle -\left\langle x,e\right\rangle
\left\langle e,y\right\rangle \right\vert \leq \frac{1}{4}\left\vert \Phi
-\phi \right\vert \left\vert \Gamma -\gamma \right\vert .  \label{1.2}
\end{equation}%
The constant $\frac{1}{4}$ is best possible in the sense that it cannot be
replaced by a smaller constant.
\end{theorem}

In \cite{SSD2}, the following refinement of (\ref{1.2}) has been pointed out.

\begin{theorem}
\label{t2}Let $H,$ $\mathbb{K}$ and $e$ be as in Theorem \ref{t1}. If $\phi
,\Phi ,\gamma ,\Gamma ,x,y$ satisfy (\ref{1.1}) or, equivalently (see \cite[%
Lemma 1]{SSD2})%
\begin{equation}
\left\Vert x-\frac{\phi +\Phi }{2}e\right\Vert \leq \frac{1}{2}\left\vert
\Phi -\phi \right\vert ,\ \ \ \left\Vert y-\frac{\gamma +\Gamma }{2}%
e\right\Vert \leq \frac{1}{2}\left\vert \Gamma -\gamma \right\vert ,
\label{1.3}
\end{equation}%
then%
\begin{multline}
\left\vert \left\langle x,y\right\rangle -\left\langle x,e\right\rangle
\left\langle e,y\right\rangle \right\vert  \label{1.4} \\
\leq \frac{1}{4}\left\vert \Phi -\phi \right\vert \left\vert \Gamma -\gamma
\right\vert -\left[ \func{Re}\left\langle \Phi e-x,x-\phi e\right\rangle %
\right] ^{\frac{1}{2}}\left[ \func{Re}\left\langle \Gamma e-y,y-\gamma
e\right\rangle \right] ^{\frac{1}{2}}.
\end{multline}
\end{theorem}

In \cite{NU}, N. Ujevi\'{c} has generalised Theorem \ref{t1} for the case of
real inner product spaces as follows.

\begin{theorem}
\label{t3}Let $\left( H,\left\langle \cdot ,\cdot \right\rangle \right) $ be
an inner product space over the real numbers field $\mathbb{R}$, and $%
\left\{ e_{i}\right\} _{i\in \left\{ 1,\dots ,n\right\} }$ an orthornormal
family in $H,$ i.e., we recall that $\left\langle e_{i},e_{j}\right\rangle =0
$ if $i\neq j$ and $\left\Vert e_{i}\right\Vert =1,$ $i,j\in \left\{ 1,\dots
,n\right\} .$ If $\phi _{i},\gamma _{i},\Phi _{i},\Gamma _{i}\in \mathbb{R}$%
, $i\in \left\{ 1,\dots ,n\right\} $ satisfy the condition%
\begin{equation}
\left\langle \sum_{i=1}^{n}\Phi _{i}e_{i}-x,x-\sum_{i=1}^{n}\phi
_{i}e_{i}\right\rangle \geq 0,\ \ \ \ \ \ \left\langle \sum_{i=1}^{n}\Gamma
_{i}e_{i}-y,y-\sum_{i=1}^{n}\gamma _{i}e_{i}\right\rangle \geq 0,
\label{1.5}
\end{equation}%
then one has the inequality:%
\begin{equation}
\left\vert \left\langle x,y\right\rangle -\sum_{i=1}^{n}\left\langle
x,e_{i}\right\rangle \left\langle e_{i},y\right\rangle \right\vert \leq 
\frac{1}{4}\left[ \sum_{i=1}^{n}\left( \Phi _{i}-\phi _{i}\right) ^{2}\cdot
\sum_{i=1}^{n}\left( \Gamma _{i}-\gamma _{i}\right) ^{2}\right] ^{\frac{1}{2}%
}.  \label{1.6}
\end{equation}%
The constant $\frac{1}{4}$ is best possible in the sense that it cannot be
replaced by a smaller constant.
\end{theorem}

We note that the key point in his proof is the following identity:%
\begin{multline}
\sum_{i=1}^{n}\left( \left\langle x,e_{i}\right\rangle -\phi _{i}\right)
\left( \Phi _{i}-\left\langle x,e_{i}\right\rangle \right) -\left\langle
x-\sum_{i=1}^{n}\phi _{i}e_{i},\sum_{i=1}^{n}\Phi _{i}e_{i}-x\right\rangle 
\label{1.7} \\
=\left\Vert x\right\Vert ^{2}-\sum_{i=1}^{n}\left\langle
x,e_{i}\right\rangle ^{2},
\end{multline}%
holding for $x\in H,$ $\phi _{i},\Phi _{i}\in \mathbb{R}$, $i\in \left\{
1,\dots ,n\right\} $ and $\left\{ e_{i}\right\} _{i\in \left\{ 1,\dots
,n\right\} }$ an orthornormal family of vectors in the real inner product
space $H.$

In this paper we point out a counterpart of Bessel's inequality in both real
and complex inner product spaces. This result will then be employed to
provide a refinement of the Gr\"{u}ss type inequality $\left( \ref{1.6}%
\right) $ for real or complex inner products. Related results as well as
integral inequalities for general measure spaces are also given.

\section{A Counterpart of Bessel's Inequality\label{s2}}

Let $\left( H,\left\langle \cdot ,\cdot \right\rangle \right) $ be an inner
product over the real or complex number field $\mathbb{K}$ and $\left\{
e_{i}\right\} _{i\in I}$ a finite or infinite family of \textit{orthornormal
vectors }in $H,$ i.e.,%
\begin{equation}
\left\langle e_{i},e_{j}\right\rangle =\left\{ 
\begin{array}{cc}
0 & \text{if \ }i\neq j \\ 
&  \\ 
1 & \text{if \ }i=j%
\end{array}%
,\right. \ \ \ i,j\in I.  \label{2.1}
\end{equation}%
It is well known that, the following inequality due to Bessel holds 
\begin{equation}
\sum_{i\in I}\left\vert \left\langle x,e_{i}\right\rangle \right\vert
^{2}\leq \left\Vert x\right\Vert ^{2}\text{ \ for any }x\in H,  \label{2.2}
\end{equation}%
where the meaning of the sum is:%
\begin{equation}
\sum_{i\in I}\left\vert \left\langle x,e_{i}\right\rangle \right\vert
^{2}:=\sup_{F\subset I}\left\{ \sum_{i\in F}\left\vert \left\langle
x,e_{i}\right\rangle \right\vert ^{2},\ \ \ F\text{ is a finite part of }%
I\right\} .  \label{2.3}
\end{equation}

The following lemma holds.

\begin{lemma}
\label{l2.1}Let $\left\{ e_{i}\right\} _{i\in I}$ be a family of
orthornormal vectors in $H,$ $F$ a finite part of $I$ and $\phi _{i},\Phi
_{i}$ $\left( i\in F\right) ,$ real or complex numbers. The following
statements are equivalent for $x\in H$

\begin{enumerate}
\item[(i)] $\func{Re}\left\langle \sum_{i\in F}\Phi _{i}e_{i}-x,x-\sum_{i\in
F}\phi _{i}e_{i}\right\rangle \geq 0$

\item[(ii)] $\left\Vert x-\sum_{i\in F}\frac{\phi _{i}+\Phi _{i}}{2}%
e_{i}\right\Vert \leq \frac{1}{2}\left( \sum_{i\in F}\left\vert \Phi
_{i}-\phi _{i}\right\vert ^{2}\right) ^{\frac{1}{2}}.$
\end{enumerate}
\end{lemma}

\begin{proof}
It is easy to see that for $y,a,A\in H,$ the following are equivalent (see 
\cite[Lemma 1]{SSD2})

\begin{enumerate}
\item[(a)] $\func{Re}\left\langle A-y,y-a\right\rangle \geq 0$ and

\item[(aa)] $\left\Vert y-\frac{a+A}{2}\right\Vert \leq \frac{1}{2}%
\left\Vert A-a\right\Vert .$
\end{enumerate}

Now, for $a=\sum_{i\in F}\phi _{i}e_{i},$ $A=\sum_{i\in F}\Phi _{i}e_{i},$
we have%
\begin{align*}
\left\Vert A-a\right\Vert & =\left\Vert \sum_{i\in F}\left( \Phi _{i}-\phi
_{i}\right) e_{i}\right\Vert =\left( \left\Vert \sum_{i\in F}\left( \Phi
_{i}-\phi _{i}\right) e_{i}\right\Vert ^{2}\right) ^{\frac{1}{2}} \\
& =\left( \sum_{i\in F}\left\vert \Phi _{i}-\phi _{i}\right\vert
^{2}\left\Vert e_{i}\right\Vert ^{2}\right) ^{\frac{1}{2}}=\left( \sum_{i\in
F}\left\vert \Phi _{i}-\phi _{i}\right\vert ^{2}\right) ^{\frac{1}{2}},
\end{align*}%
giving, for $y=x,$ the desired equivalence.
\end{proof}

The following counterpart of Bessel's inequality holds.

\begin{theorem}
\label{t2.2}Let $\left\{ e_{i}\right\} _{i\in I},$ $F,$ $\phi _{i},\Phi
_{i}, $ $i\in F$ and $x\in H$ so that either (i) or (ii) of Lemma \ref{l2.1}
holds. Then we have the inequality:%
\begin{align}
0& \leq \left\Vert x\right\Vert ^{2}-\sum_{i\in F}\left\vert \left\langle
x,e_{i}\right\rangle \right\vert ^{2}  \label{2.4} \\
& \leq \frac{1}{4}\sum_{i\in F}\left\vert \Phi _{i}-\phi _{i}\right\vert
^{2}-\func{Re}\left\langle \sum_{i\in F}\Phi _{i}e_{i}-x,x-\sum_{i\in F}\phi
_{i}e_{i}\right\rangle  \notag \\
& \leq \frac{1}{4}\sum_{i\in F}\left\vert \Phi _{i}-\phi _{i}\right\vert
^{2}.  \notag
\end{align}%
The constant $\frac{1}{4}$ is best in both inequalities.
\end{theorem}

\begin{proof}
Define%
\begin{equation*}
I_{1}:=\sum_{i\in H}\func{Re}\left[ \left( \Phi _{i}-\left\langle
x,e_{i}\right\rangle \right) \left( \overline{\left\langle
x,e_{i}\right\rangle }-\overline{\phi _{i}}\right) \right] 
\end{equation*}%
and%
\begin{equation*}
I_{2}:=\func{Re}\left[ \left\langle \sum_{i\in H}\Phi
_{i}e_{i}-x,x-\sum_{i\in H}\phi _{i}e_{i}\right\rangle \right] .
\end{equation*}%
Observe that%
\begin{equation*}
I_{1}=\sum_{i\in H}\func{Re}\left[ \Phi _{i}\overline{\left\langle
x,e_{i}\right\rangle }\right] +\sum_{i\in H}\func{Re}\left[ \overline{\phi
_{i}}\left\langle x,e_{i}\right\rangle \right] -\sum_{i\in H}\func{Re}\left[
\Phi _{i}\overline{\phi _{i}}\right] -\sum_{i\in H}\left\vert \left\langle
x,e_{i}\right\rangle \right\vert ^{2}
\end{equation*}%
and%
\begin{align*}
I_{2}& =\func{Re}\left[ \sum_{i\in H}\Phi _{i}\overline{\left\langle
x,e_{i}\right\rangle }+\sum_{i\in H}\overline{\phi _{i}}\left\langle
x,e_{i}\right\rangle -\left\Vert x\right\Vert ^{2}-\sum_{i\in H}\sum_{j\in
H}\Phi _{i}\overline{\phi _{i}}\left\langle e_{i},e_{j}\right\rangle \right] 
\\
& =\sum_{i\in H}\func{Re}\left[ \Phi _{i}\overline{\left\langle
x,e_{i}\right\rangle }\right] +\sum_{i\in H}\func{Re}\left[ \overline{\phi
_{i}}\left\langle x,e_{i}\right\rangle \right] -\left\Vert x\right\Vert
^{2}-\sum_{i\in H}\func{Re}\left[ \Phi _{i}\overline{\phi _{i}}\right] .
\end{align*}%
Consequently, subtracting $I_{2}$ from $I_{1},$ we deduce the following
equality that is useful in its turn%
\begin{multline}
\left\Vert x\right\Vert ^{2}-\sum_{i\in F}\left\vert \left\langle
x,e_{i}\right\rangle \right\vert ^{2}=\sum_{i\in H}\func{Re}\left[ \left(
\Phi _{i}-\left\langle x,e_{i}\right\rangle \right) \left( \overline{%
\left\langle x,e_{i}\right\rangle }-\overline{\phi _{i}}\right) \right] 
\label{2.5} \\
-\func{Re}\left[ \left\langle \sum_{i\in H}\Phi _{i}e_{i}-x,x-\sum_{i\in
H}\phi _{i}e_{i}\right\rangle \right] .
\end{multline}%
Using the following elementary inequality for complex numbers%
\begin{equation*}
\func{Re}\left[ a\overline{b}\right] \leq \frac{1}{4}\left\vert
a+b\right\vert ^{2},\ \ \ \ a,b\in \mathbb{K},
\end{equation*}%
for the choices $a=\Phi _{i}-\left\langle x,e_{i}\right\rangle ,$ $%
b=\left\langle x,e_{i}\right\rangle -\phi _{i}$ \ $\left( i\in F\right) ,$
we deduce%
\begin{equation}
\sum_{i\in H}\func{Re}\left[ \left( \Phi _{i}-\left\langle
x,e_{i}\right\rangle \right) \left( \overline{\left\langle
x,e_{i}\right\rangle }-\overline{\phi _{i}}\right) \right] \leq \frac{1}{4}%
\sum_{i\in H}\left\vert \Phi _{i}-\phi _{i}\right\vert ^{2}.  \label{2.6}
\end{equation}%
Making use of (\ref{2.5}), (\ref{2.6}) and the assumption (i), we deduce (%
\ref{2.4}).

The sharpness of the constant $\frac{1}{4}$ was proved for a single element $%
e,$ $\left\Vert e\right\Vert =1$ in \cite{SSD1}, or for the real case in 
\cite{NU}.

We can give here a simple proof as follows.

Assume that there is a $c>0$ such that%
\begin{align}
0& \leq \left\Vert x\right\Vert ^{2}-\sum_{i\in F}\left\vert \left\langle
x,e_{i}\right\rangle \right\vert ^{2}  \label{2.7} \\
& \leq c\sum_{i\in F}\left\vert \Phi _{i}-\phi _{i}\right\vert ^{2}-\func{Re}%
\left\langle \sum_{i\in F}\Phi _{i}e_{i}-x,x-\sum_{i\in F}\phi
_{i}e_{i}\right\rangle ,  \notag
\end{align}%
provided $\phi _{i},\Phi _{i},$ $x$ and $F$ satisfy (i) or (ii).

We choose $F=\left\{ 1\right\} ,$ $e_{1}=e_{2}=\left( \frac{1}{\sqrt{2}},%
\frac{1}{\sqrt{2}}\right) \in \mathbb{R}^{2},$ $x=\left( x_{1},x_{2}\right)
\in \mathbb{R}^{2},\ \Phi _{1}=\Phi =m>0,$ $\phi _{1}=\phi =-m,$ $H=\mathbb{R%
}^{2}$ to get from (\ref{2.7}) that%
\begin{align}
0& \leq x_{1}^{2}+x_{2}^{2}-\frac{\left( x_{1}+x_{2}\right) ^{2}}{2}
\label{2.8} \\
& \leq 4cm^{2}-\left( \frac{m}{\sqrt{2}}-x_{1}\right) \left( x_{1}+\frac{m}{%
\sqrt{2}}\right) -\left( \frac{m}{\sqrt{2}}-x_{2}\right) \left( x_{2}+\frac{m%
}{\sqrt{2}}\right) ,  \notag
\end{align}%
provided%
\begin{align}
0& \leq \left\langle me-x,x+me\right\rangle   \label{2.9} \\
& =\left( \frac{m}{\sqrt{2}}-x_{1}\right) \left( x_{1}+\frac{m}{\sqrt{2}}%
\right) +\left( \frac{m}{\sqrt{2}}-x_{2}\right) \left( x_{2}+\frac{m}{\sqrt{2%
}}\right) .  \notag
\end{align}%
If we choose $x_{1}=\frac{m}{\sqrt{2}},$ $x_{2}=-\frac{m}{\sqrt{2}},$ then (%
\ref{2.9}) is fulfilled and by (\ref{2.8}) we get $m^{2}\leq 4cm^{2},$
giving $c\geq \frac{1}{4}.$
\end{proof}

\section{A Refinement of the Gr\"{u}ss Inequality\label{s3}}

The following result holds.

\begin{theorem}
\label{t3.1}Let $\left\{ e_{i}\right\} _{i\in I}$ be a family of
orthornormal vectors in $H,$ $F$ a finite part of $I$ and $\phi _{i},\Phi
_{i},$ $\gamma _{i},\Gamma _{i}\in \mathbb{K},\ i\in F$ and $x,y\in H.$ If
either%
\begin{align}
\func{Re}\left\langle \sum_{i\in F}\Phi _{i}e_{i}-x,x-\sum_{i\in F}\phi
_{i}e_{i}\right\rangle & \geq 0,  \label{3.1} \\
\func{Re}\left\langle \sum_{i\in F}\Gamma _{i}e_{i}-y,y-\sum_{i\in F}\gamma
_{i}e_{i}\right\rangle & \geq 0,  \notag
\end{align}%
or, equivalently,%
\begin{align}
\left\Vert x-\sum_{i\in F}\frac{\Phi _{i}+\phi _{i}}{2}e_{i}\right\Vert &
\leq \frac{1}{2}\left( \sum_{i\in F}\left\vert \Phi _{i}-\phi
_{i}\right\vert ^{2}\right) ^{\frac{1}{2}},  \label{3.2} \\
\left\Vert y-\sum_{i\in F}\frac{\Gamma _{i}+\gamma _{i}}{2}e_{i}\right\Vert
& \leq \frac{1}{2}\left( \sum_{i\in F}\left\vert \Gamma _{i}-\gamma
_{i}\right\vert ^{2}\right) ^{\frac{1}{2}},  \notag
\end{align}%
hold, then we have the inequalities%
\begin{eqnarray}
&&\left\vert \left\langle x,y\right\rangle -\sum_{i\in F}\left\langle
x,e_{i}\right\rangle \left\langle e_{i},y\right\rangle \right\vert
\label{3.3} \\
&\leq &\frac{1}{4}\left( \sum_{i\in F}\left\vert \Phi _{i}-\phi
_{i}\right\vert ^{2}\right) ^{\frac{1}{2}}\cdot \left( \sum_{i\in
F}\left\vert \Gamma _{i}-\gamma _{i}\right\vert ^{2}\right) ^{\frac{1}{2}} 
\notag \\
&&-\left[ \func{Re}\left\langle \sum_{i\in F}\Phi _{i}e_{i}-x,x-\sum_{i\in
F}\phi _{i}e_{i}\right\rangle \right] ^{\frac{1}{2}}\left[ \func{Re}%
\left\langle \sum_{i\in F}\Gamma _{i}e_{i}-y,y-\sum_{i\in F}\gamma
_{i}e_{i}\right\rangle \right] ^{\frac{1}{2}}  \notag \\
&\leq &\frac{1}{4}\left( \sum_{i\in F}\left\vert \Phi _{i}-\phi
_{i}\right\vert ^{2}\right) ^{\frac{1}{2}}\cdot \left( \sum_{i\in
F}\left\vert \Gamma _{i}-\gamma _{i}\right\vert ^{2}\right) ^{\frac{1}{2}}. 
\notag
\end{eqnarray}%
The constant $\frac{1}{4}$ is best possible.
\end{theorem}

\begin{proof}
Using Schwartz's inequality in the inner product space $\left(
H,\left\langle \cdot ,\cdot \right\rangle \right) $ one has%
\begin{multline}
\left\vert \left\langle x-\sum_{i\in F}\left\langle x,e_{i}\right\rangle
e_{i},y-\sum_{i\in F}\left\langle y,e_{i}\right\rangle e_{i}\right\rangle
\right\vert ^{2}  \label{3.4} \\
\leq \left\Vert x-\sum_{i\in F}\left\langle x,e_{i}\right\rangle
e_{i}\right\Vert ^{2}\left\Vert y-\sum_{i\in F}\left\langle
y,e_{i}\right\rangle e_{i}\right\Vert ^{2}
\end{multline}%
and since a simple calculation shows that 
\begin{equation*}
\left\langle x-\sum_{i\in F}\left\langle x,e_{i}\right\rangle
e_{i},y-\sum_{i\in F}\left\langle y,e_{i}\right\rangle e_{i}\right\rangle
=\left\langle x,y\right\rangle -\sum_{i\in F}\left\langle
x,e_{i}\right\rangle \left\langle e_{i},y\right\rangle 
\end{equation*}%
and 
\begin{equation*}
\left\Vert x-\sum_{i\in F}\left\langle x,e_{i}\right\rangle e_{i}\right\Vert
^{2}=\left\Vert x\right\Vert ^{2}-\sum_{i\in F}\left\vert \left\langle
x,e_{i}\right\rangle \right\vert ^{2}
\end{equation*}%
for any $x,y\in H,$ then by (\ref{3.4}) and by the counterpart of Bessel's
inequality in Theorem \ref{t2.2}, we have%
\begin{eqnarray}
&&\left\vert \left\langle x,y\right\rangle -\sum_{i\in F}\left\langle
x,e_{i}\right\rangle \left\langle e_{i},y\right\rangle \right\vert ^{2}
\label{3.5} \\
&\leq &\left( \left\Vert x\right\Vert ^{2}-\sum_{i\in F}\left\vert
\left\langle x,e_{i}\right\rangle \right\vert ^{2}\right) \left( \left\Vert
y\right\Vert ^{2}-\sum_{i\in F}\left\vert \left\langle y,e_{i}\right\rangle
\right\vert ^{2}\right)   \notag \\
&\leq &\left[ \frac{1}{4}\sum_{i\in F}\left\vert \Phi _{i}-\phi
_{i}\right\vert ^{2}-\func{Re}\left\langle \sum_{i\in F}\Phi
_{i}e_{i}-x,x-\sum_{i\in F}\phi _{i}e_{i}\right\rangle \right]   \notag \\
&&\times \left[ \frac{1}{4}\sum_{i\in F}\left\vert \Gamma _{i}-\gamma
_{i}\right\vert ^{2}-\func{Re}\left\langle \sum_{i\in F}\Gamma
_{i}e_{i}-y,y-\sum_{i\in F}\gamma _{i}e_{i}\right\rangle \right]   \notag \\
&\leq &\left[ \frac{1}{4}\left( \sum_{i\in F}\left\vert \Phi _{i}-\phi
_{i}\right\vert ^{2}\right) ^{\frac{1}{2}}\cdot \left( \sum_{i\in
F}\left\vert \Gamma _{i}-\gamma _{i}\right\vert ^{2}\right) ^{\frac{1}{2}%
}\right.   \notag \\
&&-\left. \left[ \func{Re}\left\langle \sum_{i\in F}\Phi
_{i}e_{i}-x,x-\sum_{i\in F}\phi _{i}e_{i}\right\rangle \right] ^{\frac{1}{2}}%
\left[ \func{Re}\left\langle \sum_{i\in F}\Gamma _{i}e_{i}-y,y-\sum_{i\in
F}\gamma _{i}e_{i}\right\rangle \right] ^{\frac{1}{2}}\right]   \notag
\end{eqnarray}%
where, for the last inequality, we have made use of the inequality%
\begin{equation*}
\left( m^{2}-n^{2}\right) \left( p^{2}-q^{2}\right) \leq \left( mp-nq\right)
^{2},
\end{equation*}%
where $m,n,p,q>0.$

Taking the square root in (\ref{3.5}) and observing that the quantity in the
last square bracket is nonnegative (see for example (\ref{2.4})), we deduce
the desired result (\ref{3.3}).

The best constant has been proved in \cite{SSD1} for one element and we omit
the details.
\end{proof}

\section{Some Companion Inequalities\label{s4}}

The following companion of the Gr\"{u}ss inequality also holds.

\begin{theorem}
\label{t4.1}Let $\left\{ e_{i}\right\} _{i\in I}$ be a family of
orthornormal vectors in $H,$ $F$ a finite part of $I$ and $\phi _{i},\Phi
_{i}\in \mathbb{K},\ i\in F$ and $x,y\in H$ such that%
\begin{equation}
\func{Re}\left\langle \sum_{i\in F}\Phi _{i}e_{i}-\frac{x+y}{2},\frac{x+y}{2}%
-\sum_{i\in F}\phi _{i}e_{i}\right\rangle \geq 0  \label{4.1}
\end{equation}%
or, equivalently,%
\begin{equation}
\left\Vert \frac{x+y}{2}-\sum_{i\in F}\frac{\Phi _{i}+\phi _{i}}{2}\cdot
e_{i}\right\Vert \leq \frac{1}{2}\left( \sum_{i\in F}\left\vert \Phi
_{i}-\phi _{i}\right\vert ^{2}\right) ^{\frac{1}{2}},  \label{4.2}
\end{equation}%
then we have the inequality%
\begin{equation}
\func{Re}\left[ \left\langle x,y\right\rangle -\sum_{i\in F}\left\langle
x,e_{i}\right\rangle \left\langle e_{i},y\right\rangle \right] \leq \frac{1}{%
4}\sum_{i\in F}\left\vert \Phi _{i}-\phi _{i}\right\vert ^{2}.  \label{4.3}
\end{equation}%
The constant $\frac{1}{4}$ is best possible.
\end{theorem}

\begin{proof}
Start with the well known inequality 
\begin{equation}
\func{Re}\left\langle z,u\right\rangle \leq \frac{1}{4}\left\Vert
z+u\right\Vert ^{2},\ \ \ z,u\in H.  \label{4.4}
\end{equation}%
Since%
\begin{equation*}
\left\langle x,y\right\rangle -\sum_{i\in F}\left\langle
x,e_{i}\right\rangle \left\langle e_{i},y\right\rangle =\left\langle
x-\sum_{i\in F}\left\langle x,e_{i}\right\rangle e_{i},y-\sum_{i\in
F}\left\langle y,e_{i}\right\rangle e_{i}\right\rangle ,
\end{equation*}%
for any \thinspace $x,y\in H,$ then, by (\ref{4.4}), we get%
\begin{align}
\func{Re}\left[ \left\langle x,y\right\rangle -\sum_{i\in F}\left\langle
x,e_{i}\right\rangle \left\langle e_{i},y\right\rangle \right] & =\func{Re}%
\left[ \left\langle x-\sum_{i\in F}\left\langle x,e_{i}\right\rangle
e_{i},y-\sum_{i\in F}\left\langle y,e_{i}\right\rangle e_{i}\right\rangle %
\right]  \label{4.5} \\
& \leq \frac{1}{4}\left\Vert x-\sum_{i\in F}\left\langle
x,e_{i}\right\rangle e_{i}+y-\sum_{i\in F}\left\langle y,e_{i}\right\rangle
e_{i}\right\Vert ^{2}  \notag \\
& =\left\Vert \frac{x+y}{2}-\sum_{i\in F}\left\langle \frac{x+y}{2}%
,e_{i}\right\rangle e_{i}\right\Vert ^{2}  \notag \\
& =\left\Vert \frac{x+y}{2}\right\Vert ^{2}-\sum_{i\in F}\left\vert
\left\langle \frac{x+y}{2},e_{i}\right\rangle \right\vert ^{2}.  \notag
\end{align}%
If we apply the counterpart of Bessel's inequality in Theorem \ref{t2.2} for 
$\frac{x+y}{2},$ we may state that%
\begin{equation}
\left\Vert \frac{x+y}{2}\right\Vert ^{2}-\sum_{i\in F}\left\vert
\left\langle \frac{x+y}{2},e_{i}\right\rangle \right\vert ^{2}\leq \frac{1}{4%
}\left( \sum_{i\in F}\left\vert \Phi _{i}-\phi _{i}\right\vert ^{2}\right) ^{%
\frac{1}{2}}.  \label{4.6}
\end{equation}%
Now, by making use of (\ref{4.5}) and (\ref{4.6}), we deduce (\ref{4.3}).

The fact that $\frac{1}{4}$ is the best constant in (\ref{4.3}) follows by
the fact that if in (\ref{4.1}) we choose $x=y,$ then it becomes (i) of
Lemma \ref{l2.1}, implying (\ref{2.4}), for which, we have shown that $\frac{%
1}{4}$ was the best constant.
\end{proof}

The following corollary may be of interest if we wish to evaluate the
absolute value of%
\begin{equation*}
\func{Re}\left[ \left\langle x,y\right\rangle -\sum_{i\in F}\left\langle
x,e_{i}\right\rangle \left\langle e_{i},y\right\rangle \right] .
\end{equation*}

\begin{corollary}
\label{c4.2}With the assumptions of Theorem \ref{t4.1} and if 
\begin{equation}
\func{Re}\left\langle \sum_{i\in F}\Phi _{i}e_{i}-\frac{x\pm y}{2},\frac{%
x\pm y}{2}-\sum_{i\in F}\phi _{i}e_{i}\right\rangle \geq 0  \label{4.7}
\end{equation}%
or, equivalently%
\begin{equation}
\left\Vert \frac{x\pm y}{2}-\sum_{i\in F}\frac{\Phi _{i}+\phi _{i}}{2}\cdot
e_{i}\right\Vert \leq \frac{1}{2}\left( \sum_{i\in F}\left\vert \Phi
_{i}-\phi _{i}\right\vert ^{2}\right) ^{\frac{1}{2}},  \label{4.8}
\end{equation}%
then we have the inequality%
\begin{equation}
\left\vert \func{Re}\left[ \left\langle x,y\right\rangle -\sum_{i\in
F}\left\langle x,e_{i}\right\rangle \left\langle e_{i},y\right\rangle \right]
\right\vert \leq \frac{1}{4}\sum_{i\in F}\left\vert \Phi _{i}-\phi
_{i}\right\vert ^{2}.  \label{4.9}
\end{equation}
\end{corollary}

\begin{proof}
We only remark that, if%
\begin{equation*}
\func{Re}\left\langle \sum_{i\in F}\Phi _{i}e_{i}-\frac{x-y}{2},\frac{x-y}{2}%
-\sum_{i\in F}\phi _{i}e_{i}\right\rangle \geq 0
\end{equation*}%
holds, then by Theorem \ref{t4.1} for $\left( -y\right) $ instead of $y,$ we
have%
\begin{equation*}
\func{Re}\left[ -\left\langle x,y\right\rangle +\sum_{i\in F}\left\langle
x,e_{i}\right\rangle \left\langle e_{i},y\right\rangle \right] \leq \frac{1}{%
4}\sum_{i\in F}\left\vert \Phi _{i}-\phi _{i}\right\vert ^{2},
\end{equation*}%
showing that%
\begin{equation}
\func{Re}\left[ \left\langle x,y\right\rangle -\sum_{i\in F}\left\langle
x,e_{i}\right\rangle \left\langle e_{i},y\right\rangle \right] \geq -\frac{1%
}{4}\sum_{i\in F}\left\vert \Phi _{i}-\phi _{i}\right\vert ^{2}.
\label{4.10}
\end{equation}%
Making use of (\ref{4.3}) and (\ref{4.10}), we deduce the desired inequality
(\ref{4.9}).
\end{proof}

\begin{remark}
\label{r4.3}If $H$ is a real inner product space and $m_{i},M_{i}\in \mathbb{%
R}$ with the property that%
\begin{equation}
\left\langle \sum_{i\in F}M_{i}e_{i}-\frac{x\pm y}{2},\frac{x\pm y}{2}%
-\sum_{i\in F}m_{i}e_{i}\right\rangle \geq 0  \label{4.11}
\end{equation}%
or, equivalently,%
\begin{equation}
\left\Vert \frac{x\pm y}{2}-\sum_{i\in F}\frac{M_{i}+m_{i}}{2}\cdot
e_{i}\right\Vert \leq \frac{1}{2}\left( \sum_{i\in F}\left(
M_{i}-m_{i}\right) ^{2}\right) ^{\frac{1}{2}},  \label{4.12}
\end{equation}%
then we have the Gr\"{u}ss type inequality%
\begin{equation}
\left\vert \left\langle x,y\right\rangle -\sum_{i\in F}\left\langle
x,e_{i}\right\rangle \left\langle e_{i},y\right\rangle \right\vert \leq 
\frac{1}{4}\sum_{i\in F}\left( M_{i}-m_{i}\right) ^{2}.  \label{4.13}
\end{equation}
\end{remark}

\section{Integral Inequalities\label{s5}}

Let $\left( \Omega ,\Sigma ,\mu \right) $ be a measure space consisting of a
set $\Omega ,$ a $\sigma -$algebra of parts $\Sigma $ and a countably
additive and positive measure $\mu $ on $\Sigma $ with values in $\mathbb{R}%
\cup \left\{ \infty \right\} .$ Let $\rho \geq 0$ be a $\mu -$measurable
function on $\Omega .$ Denote by $L_{\rho }^{2}\left( \Omega ,\mathbb{K}%
\right) $ the Hilbert space of all real or complex valued functions defined
on $\Omega $ and $2-\rho -$integrable on $\Omega ,$ i.e.,%
\begin{equation}
\int_{\Omega }\rho \left( s\right) \left\vert f\left( s\right) \right\vert
^{2}d\mu \left( s\right) <\infty .  \label{5.1}
\end{equation}

Consider the family $\left\{ f_{i}\right\} _{i\in I}$ of functions in $%
L_{\rho }^{2}\left( \Omega ,\mathbb{K}\right) $ with the properties that%
\begin{equation}
\int_{\Omega }\rho \left( s\right) f_{i}\left( s\right) \overline{f_{j}}%
\left( s\right) d\mu \left( s\right) =\delta _{ij},\ \ \ i,j\in I,
\label{5.2}
\end{equation}%
where $\delta _{ij}$ is $0$ if $i\neq j$ and $\delta _{ij}=1$ if $i=j.$ $%
\left\{ f_{i}\right\} _{i\in I}$ is an orthornormal family in $L_{\rho
}^{2}\left( \Omega ,\mathbb{K}\right) .$

The following proposition holds.

\begin{proposition}
\label{p5.1}Let $\left\{ f_{i}\right\} _{i\in I}$ be an orthornormal family
of functions in $L_{\rho }^{2}\left( \Omega ,\mathbb{K}\right) ,$ $F$ a
finite subset of $I,$ $\phi _{i},\Phi _{i}\in \mathbb{K}$ $\left( i\in
F\right) $ and $f\in L_{\rho }^{2}\left( \Omega ,\mathbb{K}\right) ,$ so
that either%
\begin{equation}
\int_{\Omega }\rho \left( s\right) \func{Re}\left[ \left( \sum_{i\in F}\Phi
_{i}f_{i}\left( s\right) -f\left( s\right) \right) \left( \overline{f}\left(
s\right) -\sum_{i\in F}\overline{\phi _{i}}\text{ }\overline{f_{i}}\left(
s\right) \right) \right] d\mu \left( s\right) \geq 0  \label{5.3}
\end{equation}%
or, equivalently,%
\begin{equation}
\int_{\Omega }\rho \left( s\right) \left\vert f\left( s\right) -\sum_{i\in F}%
\frac{\Phi _{i}+\phi _{i}}{2}f_{i}\left( s\right) \right\vert ^{2}d\mu
\left( s\right) \leq \frac{1}{4}\sum_{i\in F}\left\vert \Phi _{i}-\phi
_{i}\right\vert ^{2}.  \label{5.4}
\end{equation}%
Then we have the inequality%
\begin{align}
0& \leq \int_{\Omega }\rho \left( s\right) \left\vert f\left( s\right)
\right\vert ^{2}d\mu \left( s\right) -\sum_{i\in F}\left\vert \int_{\Omega
}\rho \left( s\right) f\left( s\right) \overline{f_{i}}\left( s\right) d\mu
\left( s\right) \right\vert ^{2}  \label{5.5} \\
& \leq \frac{1}{4}\sum_{i\in F}\left\vert \Phi _{i}-\phi _{i}\right\vert ^{2}
\notag \\
& \ \ \ \ \ -\int_{\Omega }\rho \left( s\right) \func{Re}\left[ \left(
\sum_{i\in F}\Phi _{i}f_{i}\left( s\right) -f\left( s\right) \right) \left( 
\overline{f}\left( s\right) -\sum_{i\in F}\overline{\phi _{i}}\text{ }%
\overline{f_{i}}\left( s\right) \right) \right] d\mu \left( s\right)   \notag
\\
& \leq \frac{1}{4}\sum_{i\in F}\left\vert \Phi _{i}-\phi _{i}\right\vert
^{2}.  \notag
\end{align}%
The constant $\frac{1}{4}$ is best possible in both inequalities.
\end{proposition}

The proof follows by Theorem \ref{t2.2} applied for the Hilbert space $%
L_{\rho }^{2}\left( \Omega ,\mathbb{K}\right) $ and the orthornormal family $%
\left\{ f_{i}\right\} _{i\in I}.$

The following Gr\"{u}ss type inequality also holds.

\begin{proposition}
\label{p5.2}Let $\left\{ f_{i}\right\} _{i\in I}$ and $F$ be as in
Proposition \ref{p5.1}. If $\phi _{i},\Phi _{i},\gamma _{i},\Gamma _{i}\in 
\mathbb{K}$ $\left( i\in F\right) $ and $f,g\in L_{\rho }^{2}\left( \Omega ,%
\mathbb{K}\right) $ so that either%
\begin{align}
\int_{\Omega }\rho \left( s\right) \func{Re}\left[ \left( \sum_{i\in F}\Phi
_{i}f_{i}\left( s\right) -f\left( s\right) \right) \left( \overline{f}\left(
s\right) -\sum_{i\in F}\overline{\phi _{i}}\text{ }\overline{f_{i}}\left(
s\right) \right) \right] d\mu \left( s\right) & \geq 0,  \label{5.6} \\
\int_{\Omega }\rho \left( s\right) \func{Re}\left[ \left( \sum_{i\in
F}\Gamma _{i}f_{i}\left( s\right) -g\left( s\right) \right) \left( \overline{%
g}\left( s\right) -\sum_{i\in F}\overline{\gamma _{i}}\text{ }\overline{f_{i}%
}\left( s\right) \right) \right] d\mu \left( s\right) & \geq 0,  \notag
\end{align}%
or, equivalently,%
\begin{align}
\int_{\Omega }\rho \left( s\right) \left\vert f\left( s\right) -\sum_{i\in F}%
\frac{\Phi _{i}+\phi _{i}}{2}f_{i}\left( s\right) \right\vert ^{2}d\mu
\left( s\right) & \leq \frac{1}{4}\sum_{i\in F}\left\vert \Phi _{i}-\phi
_{i}\right\vert ^{2},  \label{5.7} \\
\int_{\Omega }\rho \left( s\right) \left\vert g\left( s\right) -\sum_{i\in F}%
\frac{\Gamma _{i}+\gamma _{i}}{2}f_{i}\left( s\right) \right\vert ^{2}d\mu
\left( s\right) & \leq \frac{1}{4}\sum_{i\in F}\left\vert \Gamma _{i}-\gamma
_{i}\right\vert ^{2},  \notag
\end{align}%
then we have the inequalities%
\begin{multline}
\left\vert \int_{\Omega }\rho \left( s\right) f\left( s\right) g\left(
s\right) d\mu \left( s\right) \right.   \label{5.8} \\
-\left. \sum_{i\in F}\int_{\Omega }\rho \left( s\right) f\left( s\right) 
\overline{f_{i}}\left( s\right) d\mu \left( s\right) \int_{\Omega }\rho
\left( s\right) f_{i}\left( s\right) \overline{g\left( s\right) }d\mu \left(
s\right) \right\vert 
\end{multline}%
\begin{eqnarray}
&\leq &\frac{1}{4}\left( \sum_{i\in F}\left\vert \Phi _{i}-\phi
_{i}\right\vert ^{2}\right) ^{\frac{1}{2}}\left( \sum_{i\in F}\left\vert
\Gamma _{i}-\gamma _{i}\right\vert ^{2}\right) ^{\frac{1}{2}}  \notag \\
&&-\left[ \int_{\Omega }\rho \left( s\right) \func{Re}\left[ \left(
\sum_{i\in F}\Phi _{i}f_{i}\left( s\right) -f\left( s\right) \right) \left( 
\overline{f}\left( s\right) -\sum_{i\in F}\phi \overline{_{i}}\overline{f_{i}%
}\left( s\right) \right) \right] d\mu \left( s\right) \right] ^{\frac{1}{2}}
\notag \\
&&\times \left[ \int_{\Omega }\rho \left( s\right) \func{Re}\left[ \left(
\sum_{i\in F}\Gamma _{i}f_{i}\left( s\right) -g\left( s\right) \right)
\left( \overline{g}\left( s\right) -\sum_{i\in F}\overline{\gamma _{i}}%
\overline{f_{i}}\left( s\right) \right) \right] d\mu \left( s\right) \right]
^{\frac{1}{2}}  \notag \\
&\leq &\frac{1}{4}\left( \sum_{i\in F}\left\vert \Phi _{i}-\phi
_{i}\right\vert ^{2}\right) ^{\frac{1}{2}}\left( \sum_{i\in F}\left\vert
\Gamma _{i}-\gamma _{i}\right\vert ^{2}\right) ^{\frac{1}{2}}.  \notag
\end{eqnarray}%
The constant $\frac{1}{4}$ is the best possible.
\end{proposition}

The proof follows by Theorem \ref{t3.1} and we omit the details.

\begin{remark}
\label{r5.3}Similar results may be stated if we apply the inequalities in
Section \ref{s4}. We omit the details.
\end{remark}

In the case of real spaces, the following corollaries provide much simpler
sufficient conditions for the counterpart of Bessel's inequality (\ref{5.5})
or for the Gr\"{u}ss type inequality (\ref{5.8}) to hold.

\begin{corollary}
\label{c5.4}Let $\left\{ f_{i}\right\} _{i\in I}$ be an orthornormal family
of functions in the real Hilbert space $L_{\rho }^{2}\left( \Omega ,\mathbb{K%
}\right) ,$ $F$ a finite part of $I,$ $M_{i},m_{i}\in \mathbb{R}$ \ $\left(
i\in F\right) $ and $f\in L_{\rho }^{2}\left( \Omega \right) $ so that%
\begin{equation}
\sum_{i\in F}m_{i}f_{i}\left( s\right) \leq f\left( s\right) \leq \sum_{i\in
F}M_{i}f_{i}\left( s\right) \text{ \ for \ }\mu -\text{a.e. }s\in \Omega .
\label{5.9}
\end{equation}%
Then we have the inequalities%
\begin{align}
0& \leq \int_{\Omega }\rho \left( s\right) f^{2}\left( s\right) d\mu \left(
s\right) -\sum_{i\in F}\left[ \int_{\Omega }\rho \left( s\right) f\left(
s\right) f_{i}\left( s\right) d\mu \left( s\right) \right] ^{2}  \label{5.10}
\\
& \leq \frac{1}{4}\sum_{i\in F}\left( M_{i}-m_{i}\right) ^{2}  \notag \\
& \ \ \ \ \ \ -\int_{\Omega }\rho \left( s\right) \left( \sum_{i\in
F}M_{i}f_{i}\left( s\right) -f\left( s\right) \right) \left( f\left(
s\right) -\sum_{i\in F}m_{i}f_{i}\left( s\right) \right) d\mu \left( s\right)
\notag \\
& \leq \frac{1}{4}\sum_{i\in F}\left( M_{i}-m_{i}\right) ^{2}.  \notag
\end{align}%
The constant $\frac{1}{4}$ is best possible.
\end{corollary}

\begin{corollary}
\label{c5.5}Let $\left\{ f_{i}\right\} _{i\in I}$ and $F$ be as in Corollary %
\ref{c5.4}. If $M_{i},m_{i},N_{i},n_{i}\in \mathbb{R}$ $\left( i\in F\right) 
$ and $f,g\in L_{\rho }^{2}\left( \Omega \right) $ such that%
\begin{equation}
\sum_{i\in F}m_{i}f_{i}\left( s\right) \leq f\left( s\right) \leq \sum_{i\in
F}M_{i}f_{i}\left( s\right)  \label{5.11}
\end{equation}%
and%
\begin{equation*}
\sum_{i\in F}n_{i}f_{i}\left( s\right) \leq g\left( s\right) \leq \sum_{i\in
F}N_{i}f_{i}\left( s\right) \text{ \ for \ }\mu -\text{a.e. }s\in \Omega .
\end{equation*}%
Then we have the inequalities%
\begin{eqnarray}
&&\left\vert \int_{\Omega }\rho \left( s\right) f\left( s\right) g\left(
s\right) d\mu \left( s\right) \right.  \label{5.12} \\
&&-\left. \sum_{i\in F}\int_{\Omega }\rho \left( s\right) f\left( s\right)
f_{i}\left( s\right) d\mu \left( s\right) \int_{\Omega }\rho \left( s\right)
g\left( s\right) f_{i}\left( s\right) d\mu \left( s\right) \right\vert 
\notag \\
&\leq &\frac{1}{4}\left( \sum_{i\in F}\left( M_{i}-m_{i}\right) ^{2}\right)
^{\frac{1}{2}}\left( \sum_{i\in F}\left( N_{i}-n_{i}\right) ^{2}\right) ^{%
\frac{1}{2}}  \notag \\
&&-\left[ \int_{\Omega }\rho \left( s\right) \left( \sum_{i\in
F}M_{i}f_{i}\left( s\right) -f\left( s\right) \right) \left( f\left(
s\right) -\sum_{i\in F}m_{i}f_{i}\left( s\right) \right) d\mu \left(
s\right) \right] ^{\frac{1}{2}}  \notag \\
&&\times \left[ \int_{\Omega }\rho \left( s\right) \left( \sum_{i\in
F}N_{i}f_{i}\left( s\right) -g\left( s\right) \right) \left( g\left(
s\right) -\sum_{i\in F}n_{i}f_{i}\left( s\right) \right) d\mu \left(
s\right) \right] ^{\frac{1}{2}}  \notag \\
&\leq &\frac{1}{4}\left( \sum_{i\in F}\left( M_{i}-m_{i}\right) ^{2}\right)
^{\frac{1}{2}}\left( \sum_{i\in F}\left( N_{i}-n_{i}\right) ^{2}\right) ^{%
\frac{1}{2}}.  \notag
\end{eqnarray}
\end{corollary}

\end{document}